\documentclass{article}%
\usepackage{amsfonts}
\usepackage{amsmath}
\usepackage[letterpaper]{geometry}%
\usepackage{amssymb}%
\usepackage{graphicx}
\newtheorem{theorem}{Theorem}

\begin{document}

\title{Pell's equation without irrational numbers}
\author{N.\ J. Wildberger\\School of Mathematics and Statistics\\UNSW Sydney Australia}
\maketitle

\begin{abstract}
We give a simple way to generate an infinite number of solutions to Pell's
equation $x^{2}-Dy^{2}=1$, requiring only basic matrix arithmetic and no
knowledge of irrational numbers. We illustrate the method for $D=2,7$ and
$61.$ Connections to the Stern-Brocot tree and universal geometry are also discussed.

\end{abstract}

\section{Introduction}

For $D$ a positive square-free integer, Pell's equation is%
\[
x^{2}-Dy^{2}=1.
\]
This is perhaps the most important Diophantine equation. Its history goes back
to the ancient Greeks, to Archimedes' \textit{Cattle Problem} (see
\cite{Dorrie} and \cite{Vardi}), to Brahmagupta and Bhaskara, to Fermat and
Euler, and it was Lagrange who finally established the main fact:

\begin{theorem}
Pell's equation has an infinite number of integer solutions.
\end{theorem}

Proofs are found in many books and articles on number theory, see for example
\cite{Cohn}, \cite{Barbeau}, \cite{Lenstra}. Arguments based on continued
fractions for quadratic irrationals are common. The aim of this paper is to
give a self-contained, simplified proof which only requires basic knowledge of
matrix multiplication and \textit{integer arithmetic}. So you need not know
what is meant by an \textit{irrational number}.

We follow classical ideas of Lagrange and Gauss, and also insights due to
Conway \cite{Conway} and others. The use of $L$ and $R$ matrices is not new,
see for example \cite{Poorten} and \cite{Raney}, and the computation of the
fundamental automorphism of a binary quadratic form as described in
\cite[Section 6.12]{Buchmann Vollmer} also avoids irrationalities---but our
approach is simpler and requires a minimum of background knowledge.

We show how the method works explicitly in a few special cases, including the
famous historical example of $D=61,$ and mention a palindromic aspect of the
solution which also sheds light on the negative Pell equation $x^{2}%
-Dy^{2}=-1.$ While we generate an infinite number of integer solutions, we do
not address the issue of whether or not \textit{all }solutions are obtained by
this method.

At the end of the paper we briefly connect our solution to paths in the
\textit{Stern-Brocot tree}, and mention the connections with \textit{universal
geometry}.

We need just a few basic definitions. A \textbf{quadratic form} is an
expression in variables $x$ and $y$ of the form%
\begin{equation}
Q\left(  x,y\right)  =ax^{2}+2bxy+cy^{2}\label{Quadform1}%
\end{equation}
where $a,b$ and $c$ are integers. We write $Q\equiv\left(  a,b,c\right)  $ for
short. Note that we are following Gauss' convention in using an even
coefficient of $xy.$ For a (column) vector $v=\left(  x,y\right)  ^{T}$ define%
\[
Q\left(  v\right)  \equiv Q\left(  x,y\right)  =\left(  x,y\right)
\begin{pmatrix}
a & b\\
b & c
\end{pmatrix}%
\begin{pmatrix}
x\\
y
\end{pmatrix}
=v^{T}Av.
\]
The symmetric matrix
\[
A=%
\begin{pmatrix}
a & b\\
b & c
\end{pmatrix}
\]
is the \textbf{matrix} of the quadratic form $Q$. The number $ac-b^{2}$ is the
\textbf{determinant} of $Q.$

The quadratic form $Q$ is \textbf{balanced} if $a>0$ and $c<0.$ For example,
the\textbf{\ Pell quadratic form }%
\[
Q_{D}\left(  x,y\right)  \equiv x^{2}-Dy^{2}%
\]
for $D>0$ is balanced since its matrix is%
\begin{equation}
A=%
\begin{pmatrix}
1 & 0\\
0 & -D
\end{pmatrix}
\mathbf{.}\label{MatrixA}%
\end{equation}

\section{Equivalence of forms}

If $Q$ is a quadratic form with matrix $A$ and $M$ is an invertible integral
$2\times2$ matrix whose inverse also has integer entries, then we may define a
new quadratic form $Q^{\prime}$ by
\[
Q^{\prime}\left(  v\right)  \equiv Q\left(  Mv\right)  =\left(  Mv\right)
^{T}A\left(  Mv\right)  =v^{T}\left(  M^{T}AM\right)  v
\]
with matrix $A^{\prime}=M^{T}AM.$ We say that $Q$ and $Q^{\prime}$ are
\textbf{equivalent} forms, and notice that they take on the same values, since
if $Q^{\prime}\left(  v\right)  =n$ then $Q\left(  Mv\right)  =n$ and
conversely if $Q\left(  w\right)  =m$ then $Q^{\prime}\left(  M^{-1}w\right)
=m.$ Since the assumptions on $M$ imply that $\det M=\pm1,$ we must have
\[
\det A^{\prime}=\det A
\]
so equivalent forms have the same determinant.

For example, the Pell quadratic form $Q_{3}=\left(  1,0,-3\right)  $ is
equivalent to the quadratic form $Q^{\prime}=\left(  -2,-3,-3\right)  $ since
\[%
\begin{pmatrix}
1 & 0\\
1 & 1
\end{pmatrix}
^{T}%
\begin{pmatrix}
1 & 0\\
0 & -3
\end{pmatrix}%
\begin{pmatrix}
1 & 0\\
1 & 1
\end{pmatrix}
=\allowbreak%
\begin{pmatrix}
-2 & -3\\
-3 & -3
\end{pmatrix}
.
\]
The form $Q_{3}$ takes the value $\allowbreak1$ at $\left(  2,1\right)  ^{T},$
and $Q^{\prime}$ takes the same value $1$ at $\left(  2,-1\right)  ^{T}.$ The
relationship between these two solutions is
\[%
\begin{pmatrix}
1 & 0\\
1 & 1
\end{pmatrix}%
\begin{pmatrix}
2\\
-1
\end{pmatrix}
=%
\begin{pmatrix}
2\\
1
\end{pmatrix}
.
\]
Note that while $Q_{3}$ is balanced, $Q^{\prime}$ is not.

In practice we will focus on two very simple equivalences, using the following
matrices:%
\[%
\begin{tabular}
[c]{lllll}%
$L\equiv%
\begin{pmatrix}
1 & 0\\
1 & 1
\end{pmatrix}
$ &  & \textrm{and} &  & $R\equiv%
\begin{pmatrix}
1 & 1\\
0 & 1
\end{pmatrix}
.$%
\end{tabular}
\]

For a symmetric matrix
\[
A=%
\begin{pmatrix}
a & b\\
b & c
\end{pmatrix}
\]
we compute%
\begin{equation}
L^{T}AL=%
\begin{pmatrix}
1 & 0\\
1 & 1
\end{pmatrix}
^{T}%
\begin{pmatrix}
a & b\\
b & c
\end{pmatrix}%
\begin{pmatrix}
1 & 0\\
1 & 1
\end{pmatrix}
=\allowbreak%
\begin{pmatrix}
a+2b+c & b+c\\
b+c & c
\end{pmatrix}
\label{Up}%
\end{equation}
and%
\begin{equation}
R^{T}AR=%
\begin{pmatrix}
1 & 1\\
0 & 1
\end{pmatrix}
^{T}%
\begin{pmatrix}
a & b\\
b & c
\end{pmatrix}%
\begin{pmatrix}
1 & 1\\
0 & 1
\end{pmatrix}
=\allowbreak%
\begin{pmatrix}
a & a+b\\
a+b & a+2b+c
\end{pmatrix}
.\label{Down}%
\end{equation}

Going from $A$ to $L^{T}AL$ we call a\textbf{\ left step}. Going from $A$ to
$R^{T}AR$ we call a \textbf{right step}.

\section{Solving Pell's equation}

Here is our strategy for generating an infinite number of solutions to Pell's
equation $x^{2}-Dy^{2}=1$.

\begin{enumerate}
\item First observe that $e=\left(  1,0\right)  ^{T}$ is a solution.

\item Perform a \textit{suitable (non-empty) sequence of left and right
steps}, beginning with the Pell quadratic form $Q_{D}=\left(  1,0,-D\right)  $
with matrix (\ref{MatrixA}), and ending also with $Q_{D}=\left(
1,0,-D\right)  .$ This yields a matrix equation of the form%
\begin{equation}
N^{T}AN=A\label{Nequation}%
\end{equation}
for some invertible matrix $N$ with \textit{strictly positive entries},
actually of the form%
\[
N=%
\begin{pmatrix}
u & Dv\\
v & u
\end{pmatrix}
.
\]

\item Then observe that $Ne=\left(  u,v\right)  ^{T},$ the first column of
$N,$ is also a solution, since
\[
\left(  Ne\right)  ^{T}A\left(  Ne\right)  =e^{T}\left(  N^{T}AN\right)
e=e^{T}Ae=1.
\]

\item Iterate to generate an infinite family of solutions, namely
$e,Ne,N^{2}e,\cdots.$
\end{enumerate}

Let us now explain what the \textit{suitable sequence of left and right steps
is}, and why we can be sure to return to $Q_{D}$ after a finite number of such
steps. The key is simply to \textit{ensure that at each stage we obtain a
balanced form}.

If at some stage we have a balanced form $\left(  a,b,c\right)  $, we compute
the \textbf{total}
\[
T\equiv a+2b+c.
\]
If $T>0$ we \textit{perform a left step}. If $T<0$ we \textit{perform a right
step}. Then (\ref{Up}) and (\ref{Down}) ensure that the new form $\left(
a^{\prime},b^{\prime},c^{\prime}\right)  $ is also balanced, and it is clearly
equivalent to $\left(  a,b,c\right)  $.

The case $T=0\ $cannot arise if the determinant of $\left(  a,b,c\right)  $ is
$-D,$ since the identity%
\[
-D=ac-b^{2}=\left(  a+2b+c\right)  c-\left(  b+c\right)  ^{2}%
\]
would then imply that $D$ is a square.

To record the step at each stage, we write either%
\[
\left(  a,b,c\right)  L\left(  a+2b+c,b+c,c\right)
\]
or
\[
\left(  a,b,c\right)  R\left(  a,a+b,a+2b+c\right)  .
\]

Suppose we have taken a step from the balanced quadratic form $\left(
a,b,c\right)  $ to the balanced quadratic form $\left(  a^{\prime},b^{\prime
},c^{\prime}\right)  .$ Can we tell whether this step was left or right? In
the case of a left step $a^{\prime}-2b^{\prime}+c^{\prime}=a>0$ while in the
case of a right step $a^{\prime}-2b^{\prime}+c^{\prime}=c<0$. So the total of
the quadratic form $\left(  a^{\prime},-b^{\prime},c^{\prime}\right)  $
determines whether the step was left or right.

If we start with the quadratic form $Q_{D}=\left(  1,0,-D\right)  $ with
determinant $-D<0,$ and arrive after a certain number of steps to the
quadratic form $\left(  a,b,c\right)  $ then we must have $a>0$, $c<0$ and
\[
ac-b^{2}=-D.
\]
Since $ac$ is negative, there are only a finite number of integers $a,b$ and
$c$ that satisfy this equation. So we must eventually return to a form that we
have encountered before, and since we are capable of finding the inverse of
any step, we must return first to the initial form, namely $Q_{D}.$

At the beginning with $Q_{D}$ the total is $T=1-D<0$ so the first step must be
a right step. Further right steps increase the last entry of $\left(
1,b,c\right)  ,$ so before we return to $Q_{D}$ we must have taken at least
one left step. So the cumulative matrix $N,$ which is the product of the
matrices $R$ and $L$ corresponding to the steps we have taken, has strictly
positive entries. This ensures that the sequence of solutions $e,Ne,N^{2}%
e,\cdots$ is increasing in each component, and so really does represent an
infinite family of solutions.

If%
\[
N=%
\begin{pmatrix}
\alpha & \beta\\
\gamma & \delta
\end{pmatrix}
\]
then more generally (\ref{Nequation}) gives us the identity
\[
\left(  \alpha x+\beta y\right)  ^{2}-D\left(  \gamma x+\delta y\right)
^{2}=x^{2}-Dy^{2}%
\]
which allows us to generate new solutions to $Q\left(  x,y\right)  =n$ from
existing ones, for any integer $n.$ It also allows us to deduce that $N$ has
the form
\[
N=%
\begin{pmatrix}
u & Dv\\
v & u
\end{pmatrix}
.
\]

\section{Example of $D=2$}

The simplest example of Pell's equation is
\[
x^{2}-2y^{2}=1.
\]
The ancient Greeks knew how to generate solutions. Geometrically we are
looking for integer points on a hyperbola, as shown.
\begin{center}
\includegraphics[
height=4.9527cm,
width=7.7265cm
]%
{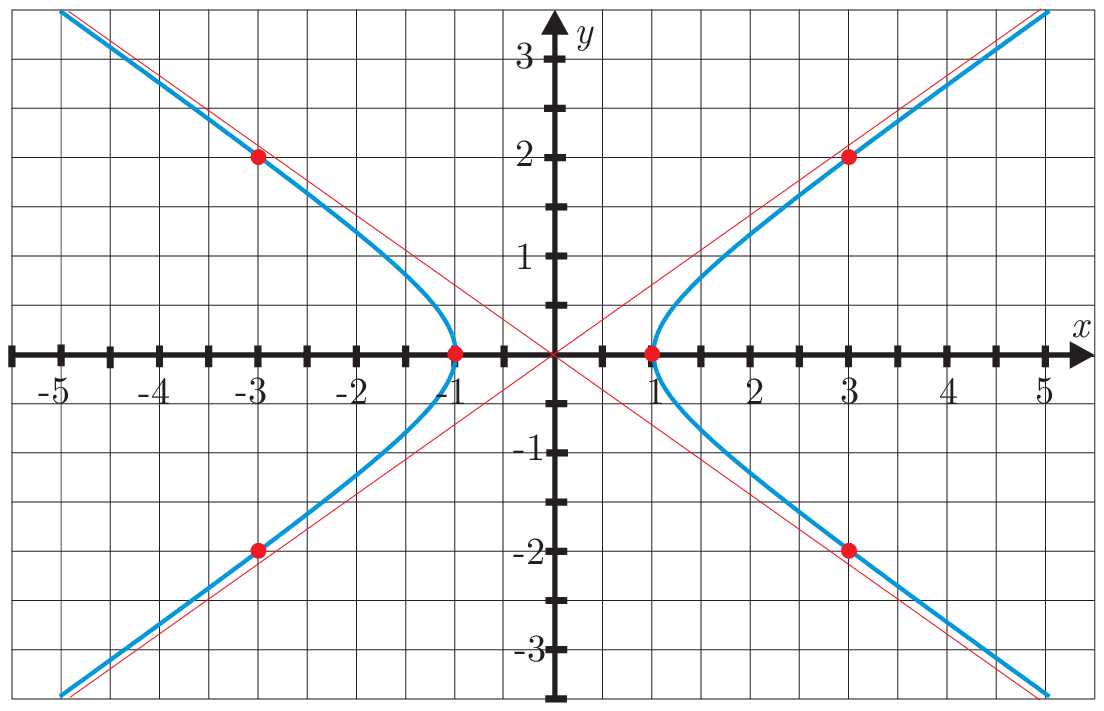}%
\\
The hyperbola $x^{2}-2y^{2}=1$%
\label{Pells2}%
\end{center}
We begin with the Pell quadratic form $Q_{2}=\left(  1,0,-2\right)  $ with a
total of $T=1+2\times0-2=-1<0$ so we take a right step to get%
\[
\left(  1,0,-2\right)  R\left(  1,1,-1\right)  .
\]
This corresponds to the matrix equation
\[
R^{T}%
\begin{pmatrix}
1 & 0\\
0 & -2
\end{pmatrix}
R=\allowbreak%
\begin{pmatrix}
1 & 1\\
1 & -1
\end{pmatrix}
.
\]
The total is now $T=1+2\times1-1=\allowbreak2>0$ so we take a left step to get%
\[
\left(  1,0,-2\right)  R\left(  1,1,-1\right)  L\left(  2,0,-1\right)  .
\]
The total is now $T=1$ so we take another left step%
\[
\left(  1,0,-2\right)  R\left(  1,1,-1\right)  L\left(  2,0,-1\right)
L\left(  1,-1,-1\right)  .
\]
The total is now $T=-2<0$ so we take a right step%
\begin{equation}
\left(  1,0,-2\right)  R\left(  1,1,-1\right)  L\left(  2,0,-1\right)
L\left(  1,-1,-1\right)  R\left(  1,0,-2\right) \label{D2}%
\end{equation}
to return to our starting quadratic form $Q_{2}.$ This yields the matrix
equation%
\[
\left(  RL^{2}R\right)  ^{T}%
\begin{pmatrix}
1 & 0\\
0 & -2
\end{pmatrix}
RL^{2}R=%
\begin{pmatrix}
1 & 0\\
0 & -2
\end{pmatrix}
.
\]
So
\[
N=RL^{2}R=\allowbreak%
\begin{pmatrix}
3 & 4\\
2 & 3
\end{pmatrix}
\]
gives the non-trivial solution $v=Ne=\left(  3,2\right)  ^{T}$ to Pell's
equation, and by taking higher powers of $N$ we get further solutions $\left(
17,12\right)  ^{T},$ $\left(  99,70\right)  ^{T},$ $\left(  577,408\right)
^{T}$ and so on.

More generally the form of $N$ establishes the identity%
\[
\left(  3x+4y\right)  ^{2}-2\left(  2x+3y\right)  ^{2}=\allowbreak
x^{2}-2y^{2}%
\]
which allows us to generate an infinite number of solutions to $Q_{2}\left(
x,y\right)  =n$ from any initial solution.

\section{Example of $D=7$}

The quadratic form $Q_{7}=\left(  1,0,-7\right)  $ gives the stepping sequence%
\begin{align*}
& \left(  1,0,-7\right)  R\left(  1,1,-6\right)  R\left(  1,2,-3\right)
L\left(  2,-1,-3\right) \\
& R\left(  2,1,-3\right)  L\left(  1,-2,-3\right)  R\left(  1,-1,-6\right)
R\left(  1,0,-7\right)  .
\end{align*}
This corresponds to the matrix equation%
\[
\left(  R^{2}LRLR^{2}\right)  ^{T}%
\begin{pmatrix}
1 & 0\\
0 & -7
\end{pmatrix}
R^{2}LRLR^{2}=%
\begin{pmatrix}
1 & 0\\
0 & -7
\end{pmatrix}
.
\]
The first column of
\[
N=R^{2}LRLR^{2}=\allowbreak%
\begin{pmatrix}
8 & 21\\
3 & 8
\end{pmatrix}
\]
gives the solution $\left(  8,3\right)  ^{T},$ and higher powers give also
$\left(  127,48\right)  ^{T},$ $\left(  2024,765\right)  ^{T}$ and so on.
Furthermore we have the identity%
\[
\left(  8x+21y\right)  ^{2}-7\left(  3x+8y\right)  ^{2}=\allowbreak
x^{2}-7y^{2}%
\]
which generates further solutions to any solution of $x^{2}-7y^{2}=n.$

\section{Palindromic aspect}

The palindromic nature of the sequences occurring for $N$ in the above two
examples is a general feature, and follows from the fact that the inverse of a
left step is%
\[
\left(  a^{\prime},b^{\prime},c^{\prime}\right)  L^{-1}\left(  a^{\prime
}-2b^{\prime}+c^{\prime},b^{\prime}-c^{\prime},c^{\prime}\right)
\]
while the inverse of a right step is
\[
\left(  a^{\prime},b^{\prime},c^{\prime}\right)  R^{-1}\left(  a^{\prime
},b^{\prime}-a^{\prime},a^{\prime}-2b^{\prime}+c^{\prime}\right)  .
\]
So if we negate the middle $b$ components, the sequence reads the same from
the end as from the beginning. To determine when we have arrived at the
half-way point we look for symmetry: two adjacent or almost adjacent
(separated by one) quadratic forms which differ only in the sign of the middle term.

The central symmetry allows us to identify when the \textit{negative} Pell
equation $x^{2}-Dy^{2}=-1$ is soluble, and give a positive solution for it, if
such exists. This will happen when the central part of the sequence looks like%
\[
\left(  a^{\prime},b^{\prime},c^{\prime}\right)  \left(  D,0,-1\right)
\left(  a^{\prime},-b^{\prime},c^{\prime}\right)
\]

In that case $N$ can be written as $N=MM^{T}$ with $M$ of the form%
\[
M=%
\begin{pmatrix}
Dv_{1} & u_{1}\\
u_{1} & v_{1}%
\end{pmatrix}
\]
so that its second column gives $\left(  u_{1},v_{1}\right)  $ satisfying
$u_{1}^{2}-Dv_{1}^{2}=-1.$

For the $D=2$ case, (\ref{D2}) shows that
\[
M=RL=\allowbreak%
\begin{pmatrix}
2 & 1\\
1 & 1
\end{pmatrix}
\]
giving $\left(  u_{1},v_{1}\right)  =\left(  1,1\right)  $. The $D=7$ case
does not have this form, and indeed $x^{2}-7y^{2}=-1$ has no solutions, as
$-1$ is not a square mod $7.$

We will further exploit this palindromic phenomenon in our last example.

\section{Example $D=61$}

The Pell's equation $x^{2}-61y^{2}=1$ was recognized by Fermat as being
particularly challenging, and he posed it as a problem to his contemporaries,
which was subsequently solved by Euler. However according to \cite{Barbeau}
the eleventh century Indian mathematician Bhaskara solved this problem. We
shorten the presentation by combining powers of $L$ and $R$ when convenient.
The required sequence of left and right steps exhibits an additional doubly
palindromic aspect, in that each half is itself palindromic in an opposite
way, that is with $L$ and $R$ symbols interchanged:%
\begin{align*}
& \left(  1,0,-61\right)  R^{7}\left(  1,7,-12\right)  L\left(
3,-5,-12\right)  R^{4}\left(  3,7,-4\right)  L^{3}\left(  9,-5,-4\right) \\
& R\left(  9,4,-5\right)  L\left(  12,-1,-5\right)  L\left(  5,-6,-5\right)
R\left(  5,-1,-12\right)  R\left(  5,4,-9\right) \\
& L\left(  4,-5,-9\right)  R^{3}\left(  4,7,-3\right)  L^{4}\left(
12,-5,-3\right)  R\left(  12,7,-1\right)  L^{7}\left(  61,0,-1\right) \\
& L^{7}\left(  12,-7,-1\right)  R\left(  12,5,-3\right)  L^{4}\left(
4,-7,-3\right)  R^{3}\left(  4,5,-9\right)  L\left(  5,-4,-9\right) \\
& R\left(  5,1,-12\right)  R\left(  5,6,-5\right)  L\left(  12,1,-5\right)
L\left(  9,-4,-5\right)  R\left(  9,5,-4\right) \\
& L^{3}\left(  3,-7,-4\right)  R^{4}\left(  3,5,-12\right)  L\left(
1,-7,-12\right)  R^{7}\left(  1,0,-61\right)
\end{align*}
Computation of the required matrix $N$ is shortened by the doubly palindromic
aspect:
\begin{align*}
N  & =R^{7}LR^{4}L^{3}RL^{2}R^{2}LR^{3}L^{4}RL^{7}L^{7}RL^{4}R^{3}LR^{2}%
L^{2}RL^{3}R^{4}LR^{7}\\
& =%
\begin{pmatrix}
453 & 164\\
58 & 21
\end{pmatrix}%
\begin{pmatrix}
453 & 58\\
164 & 21
\end{pmatrix}%
\begin{pmatrix}
21 & 58\\
164 & 453
\end{pmatrix}%
\begin{pmatrix}
21 & 164\\
58 & 453
\end{pmatrix}
\\
& =%
\begin{pmatrix}
1766\,319\,049 & 13\,\allowbreak795\,392\,780\\
226\,153\,980 & 1766\,319\,049
\end{pmatrix}
\end{align*}
and its first column gives the smallest non-trivial solution
$x=1766\,319\,049$ and $y=226\,153\,980.$ The occurrence of $\left(
61,0,-1\right)  $ in the middle of the sequence shows that the negative Pell
equation has a solution, and indeed the second column of
\[
M=%
\begin{pmatrix}
453 & 164\\
58 & 21
\end{pmatrix}%
\begin{pmatrix}
453 & 58\\
164 & 21
\end{pmatrix}
=\allowbreak%
\begin{pmatrix}
232\,105 & 29\,718\\
29\,718 & 3805
\end{pmatrix}
\]
satisfies $\left(  29\,718\right)  ^{2}-61\left(  3805\right)  ^{2}%
=\allowbreak-1.$

\section{Speeding up the algorithm}

The previous example suggests that for large $D,$ it is more efficient to work
out how many \textit{successive} $R$ or $L$ steps are needed at each stage.
Since
\[
R^{n}=%
\begin{pmatrix}
1 & n\\
0 & 1
\end{pmatrix}
\qquad\mathrm{and}\qquad L^{m}=%
\begin{pmatrix}
1 & 0\\
m & 1
\end{pmatrix}
\]
it is easy to compute that for arbitrary natural numbers $n$ and $m,$%
\begin{align*}
\left(  a,b,c\right)  R^{n}  & =\left(  a,b+an,c+2bn+an^{2}\right) \\
\left(  a,b,c\right)  L^{m}  & =\left(  a+2bm+cm^{2},b+cm,c\right)  .
\end{align*}
So if a right step $R$ is needed, namely if $T=a+2b+c<0,$ then we find the
smallest $k$ so that $c+2bk+ak^{2}>0,$ and choose $n=k-1.$ If a left step $L$
is needed, namely if $T=a+2b+c>0,$ then we find the smallest $l$ so that
$a+2bl+cl^{2}<0,$ and choose $m=l-1.$

So for example at the first stage of the $D=61$ case, $\left(  1,0,-61\right)
$ has a total of $T=-61<0$ and since $\left(  -61\right)  +k^{2}>0$ first when
$k=8,$ we set $n=8-1=7$, giving%
\[
\left(  1,0,-61\right)  R^{7}=\left(  1,7,-12\right)  .
\]
Now we know the next step must be $L^{m}$ for some $m,$ and we compute that
$1+2\times7l-12l^{2}<0$ first when $l=2,$ so that $m=2-1=1,$ giving
\[
\left(  1,0,-61\right)  R^{7}\left(  1,7,-12\right)  L=\left(
3,-5,-12\right)
\]
and so on.

\section{Connections with the Stern-Brocot tree}

The Stern-Brocot tree (see \cite{GKP} for a thorough discussion) is generated
from the two `fractions' $0/1$ and $1/0$ by repeatedly performing the mediant
operation. It naturally associates binary sequences in $L$ and $R $ with
strictly positive reduced fractions, as shown for the correspondence
\[
RL^{2}R\leftrightarrow\frac{7}{5}.
\]%
\begin{center}
\includegraphics[
height=3.9911cm,
width=10.2093cm
]%
{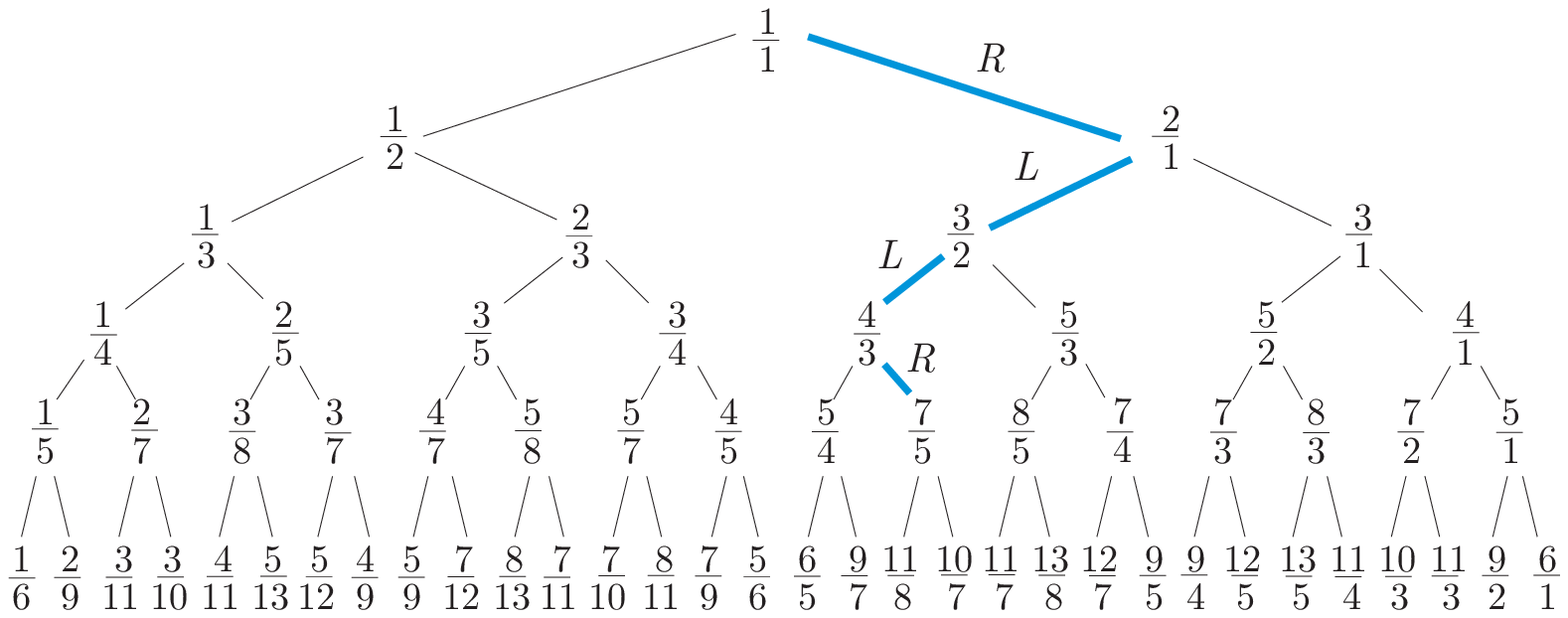}%
\\
The Stern Brocot tree
\label{SB}%
\end{center}
The Stern-Brocot tree is closely related to the Euclidean algorithm and the
theory of continued fractions, and there is a useful alternate form which
contains $2\times2$ matrices instead of fractions. The matrix corresponding to
a fraction $a/b$ has columns which are the two immediate left and right
predecessors of $a/b$ in the original tree, and whose columns sum to give the
fraction, regarded as a vector. A binary sequence of $L$'s and $R$'s can be
directly translated to a matrix by interpreting $L$ and $R$ to be the two
matrices we used earlier, and simply performing the matrix multiplication.

So computing $N$ in the Pell's equation solution amounts to following a path
in the matrix form of the Stern-Brocot tree. Each step involves adding one of
the existing columns to the other, either to the left or to the right. Shown
is the case%
\[
N=RL^{2}L=%
\begin{pmatrix}
3 & 4\\
2 & 3
\end{pmatrix}
\]
which occurred in the solution of $x^{2}-2y^{2}=1.$%
\begin{center}
\includegraphics[
height=4.2616cm,
width=11.4875cm
]%
{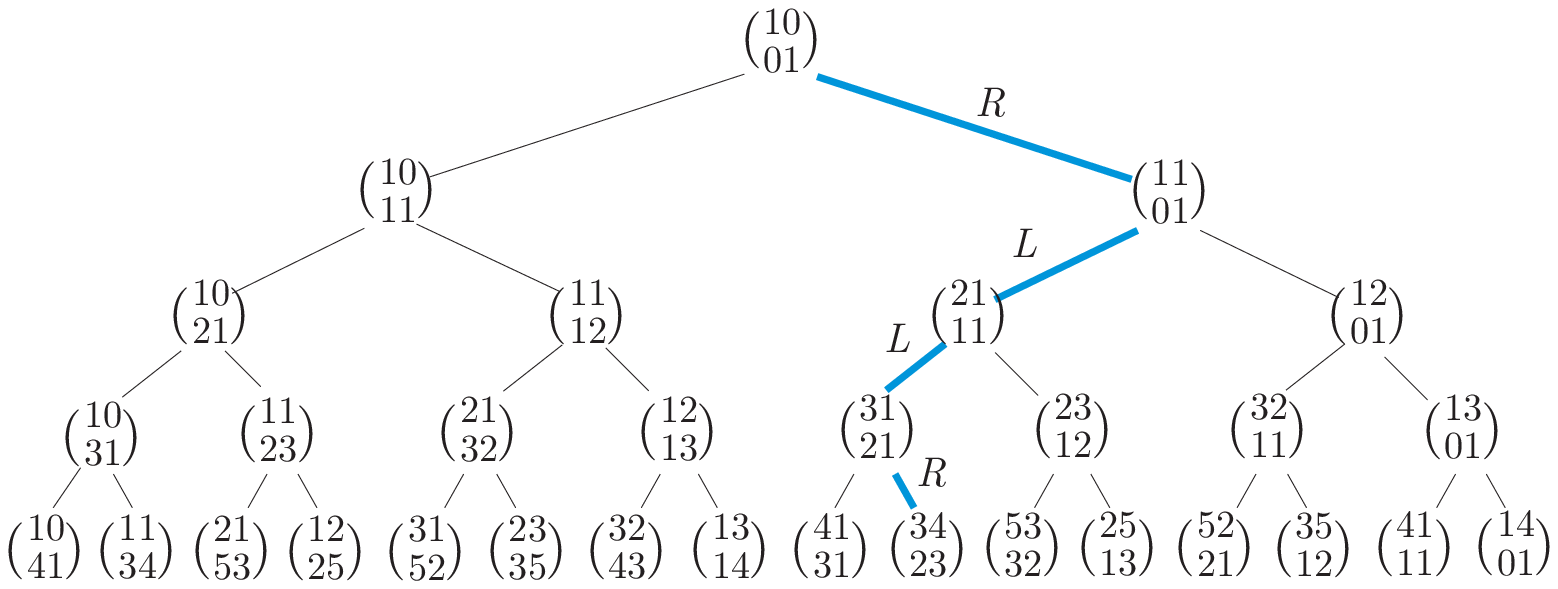}%
\\
Matrix form of the Stern-Brocot tree
\label{SB2}%
\end{center}

\section{Geometric interpretation}

We have been able to avoid irrationalities because numbers of the form
$a+b\sqrt{D}$ behave arithmetically exactly as do matrices of the form%
\[%
\begin{pmatrix}
a & Db\\
b & a
\end{pmatrix}
.
\]
This will be familiar to most readers in the imaginary case of $D=-1,$ but it
holds more generally, see for example \cite{Poorten}. It allows us to consider
a quadratic extension of the rationals not only as a subfield of the complex
numbers, but also as \textit{a field of rotation/dilation operators} on a
two-dimensional (non-Euclidean) geometry \textit{purely over the rational
numbers}.

The usual approach to metrical geometry is to consider distance and angle as
the fundamental concepts, but this is neither in the original spirit of the
ancient Greeks, nor an approach that lends itself easily to generalization.
Taking the \textit{quadrance} $x^{2}+y^{2}$ as the fundamental measurement,
instead of its square root, together with an algebraic replacement of angle
called the \textit{spread} between two lines, allows the development of a
Euclidean theory of geometry that extends both to general fields, and also to
arbitrary quadratic forms. In particular there is a very rich theory of
Euclidean geometry over the rational numbers!

This new \textit{universal geometry} is described in \cite{Wild} and
\cite{Wild2}. Our discussion of Pell's equation thus has a direct
\textit{geometric interpretation }(actually this is the context in which it
was discovered)\textit{. }So for example the matrix
\[
N=%
\begin{pmatrix}
8 & 21\\
3 & 8
\end{pmatrix}
\]
in the case of $D=7$ amounts to a \textit{rotation} in the geometry for
which\textit{\ quadrance} is defined by $Q\left(  x,y\right)  =x^{2}-7y$.

Finally we point out that a quadratic irrational number such as $\sqrt{2}$
manifests itself as a periodic infinite path down the Stern-Brocot tree. A
more general irrational number can be defined as a more general infinite path
down the Stern-Brocot tree. This approach to the definition of irrationalities
has the significant advantages that it involves no ambiguities related to
equivalences, and that the initial paths of an irrational give automatically
the best rational approximations to it.

\medskip

\textbf{2010 AMS Subject classification}: Primary---11Axx, Secondary---14Axx

\textbf{Key words and phrases}: Pell's equation, quadratic form, Stern-Brocot
tree, universal geometry


\begin{thebibliography}{99}                                                                                               %
\bibitem {Cohn}H. Cohn, \textit{Advanced Number Theory}, New York: Dover, pp.
110--111, 1980.

\bibitem {Barbeau}E. J. Barbeau, \textit{Pell's Equation}, Problem Books in
Mathematics. Springer-Verlag, 2003.

\bibitem {Buchmann Vollmer}J. Buchmann and U. Vollmer, \textit{Binary
Quadratic Forms: An Algorithmic Approach, }Algorithms and Computation in
Mathematics 20, Springer, Berlin 2007.\textit{\ }

\bibitem {Conway}J. H. Conway, \textit{The Sensual (Quadratic) Form}, Carus
Mathematical Monographs 26, Mathematical Association of America, 2005.

\bibitem {Dorrie}H. Dorrie, \textit{100 Great Problems of Elementary
Mathematics: Their history and solution}, translated (from German to English)
by David Antin, Dover, New York, 1965.

\bibitem {GKP}R. L. Graham. D. E. Knuth and O. Patashnik, \textit{Concrete
Mathematics: A Foundation for Computer Science}, Addison-Wesley, Reading, 1994.

\bibitem {Lenstra}H. W. Lenstra Jr., Solving the Pell Equation,
\textit{Notices Amer. Math. Soc. }\textbf{49}, (2002), 182--192.

\bibitem {Poorten}A. J. van der Poorten, An introduction to continued
fractions, in \textit{Diophantine Analysis}, LMS Lecture Note Ser. 109,
Cambridge University Press, Cambridge 1985, 99--138

\bibitem {Raney}S. Raney, On continued fractions and finite automata,
\textit{Math. Annalen} \textbf{206} (1973), 265--283.

\bibitem {Vardi}I. Vardi, Archimedes' Cattle Problem, \textit{Amer. Math.
Monthly }\textbf{105}, (1998), 305--319.

\bibitem {Wild}N. J. Wildberger, \textit{Divine Proportions: Rational
Trigonometry to Universal Geometry}, Wild Egg Books, Sydney, 2005.

\bibitem {Wild2}N. J. Wildberger, `One dimensional metrical geometry',
\textit{Geometriae Dedicata}, \textbf{128}, no. 1, (2007), 145--166.
\end{thebibliography}
\end{document}